\documentstyle[12pt]{article}
\begin{document} 
\setlength{\textwidth}{6in}
\setlength{\textheight}{8in}
\setlength{\topmargin}{0.3in}
\setlength{\headheight}{0in}
\setlength{\headsep}{0in}
\setlength{\parskip}{0pt}

\newcounter{probnumt}
\partopsep 0pt
\settowidth{\leftmargin}{1.}\addtolength{\leftmargin}{\labelsep}
\newenvironment{probt}{\begin{quote}\begin{enumerate} 
\setcounter{enumi}
{\value{probnumt}}}%
{  \setcounter{probnumt}
{\value{enumi}}\end{enumerate}\end{quote}}

\input{mssymb}
\newcommand{\RR}{{\Bbb R}}
\newcommand{\CC}{{\Bbb C}}
\newcommand{\NN}{{\Bbb N}}
\renewcommand{\limsup}{\overline {\lim}\,}
\renewcommand{\liminf}{\underline {\lim}\,}
\newtheorem{thm}{Theorem}
\newtheorem{lem}[thm]{Lemma}
\newtheorem{cor}[thm]{Corollary}
\newtheorem {ex}[thm]{Example}
\newtheorem{pro}{Question}
\newenvironment{rem}{\medskip\par\noindent {\bf Remark 1.}}

\newenvironment{proof}{\medskip\par\noindent{\bf Proof.}}{\hfill
$\Box$ \medskip \par}

\setcounter{page}{1}

%
%

\title{A remark of contraction semigroups on 
Banach spaces}

\author { Pei-Kee Lin\\
Department of Mathematics \\
Memphis State University\\
Memphis TN, 38152\\
USA}

\maketitle

\begin{abstract}
{Let $X$ be a complex Banach space and let $J:X \to X^*$ be a duality
section on $X$ (i.e. $\langle
x,J(x)\rangle=\|J(x)\|\|x\|=\|J(x)\|^2=\|x\|^2$).  
For any unit vector $x$ and any
($C_0$) contraction semigroup $T=\{e^{tA}:t \geq 0\}$, Goldstein
proved that if $X$ is a Hilbert space and if $|\langle T(t)
x,J(x)\rangle| \to 1 $ as $t \to \infty$, then $x$ is an
eigenvector of $A$ corresponding to a purely imaginary eigenvalue. 
In this article, we prove the similar result holds  if  $X$ is a
strictly convex complex Banach space.} \end{abstract} \newpage

Let $X$ be a complex Banach space.  Recall 
$x^* \in X^*$ is called {\em dual} of $x \in X$, if
\[\langle x, x^*\rangle=\|x^*\|\|x\|=\|x^*\|^2=\|x\|^2.\]
a mapping $J:X \to
X^*$ is said to be a {\em duality section}
if $J(x)$ is dual to $x$ for all $x \in X$.
It is known that if $X$ is a Hilbert space, then $J(x)=x$.
In \cite {Gb}, Goldstein proved the following theorem.

\begin{thm}\label{A} Let $A$ generate a $(C_0)$ contraction
semigroup $T=\{e^{tA}:t \geq 0\}$ on a complex Banach  space.
If $X$ is a Hilbert space, then for any unit vector $x \in X$,
\[\lim_{t \to \infty}|\langle T(t)x, J(x) \rangle| = 1\]
implies that $A \, x=i \lambda x$ for some real number $\lambda$.
\end{thm}

He also proved that Theorem~\ref {A} is
not true if $X$ replace by an $L_\infty$-space.  In this article, we
prove that
Theorem~\ref {A} is true if and only if $X$ is strictly convex.

Recall a Banach space is said to be {\em strictly convex} if for
any  two unit vectors $x,y \in X$, either $x=y$ or $\|x+y\| <2$.
By Hahn-Banach Theorem, $X$ is strictly convex if and only if for
any unit vectors $x,y$ and any nonzero linear functional $x^*$,
$\langle x, x^*\rangle=\|x^*\|=\langle y, x^*\rangle$ implies
$x=y$.

\begin{ex}\label {B}  Suppose $x, y$ are two unit vectors on $X$
and $x^*$ is a linear functional on $X$ such that
\[ \langle x, x^*\rangle =1=\|x^*\|=\langle y, x^*\rangle.\]
Let $A$ be defined by
\[A(z)=i \langle z, x^*\rangle x
+ ( i -\frac 1 2) (z-\langle z, x^*\rangle x).\]
Then   for any $z$,
\[e^{tA} z=e^{it}\langle z, x^*\rangle x +e^{it} e^{-\frac t 2}
(z-\langle z, x^*\rangle x)
=e^{it} \left(e^{-\frac t 2} z +(1-e^{-\frac t 2}) \langle z, x^*\rangle x
\right).\]
Hence, $A$ generates a $(C_0)$ contraction semigroup $T=\{e^{tA}:
t \geq 0\}$.  On the other hand,
\[
\langle T(t) y, x^* \rangle = e^{it}
\left \langle \langle y, x^*\rangle x + e^{-\frac t 2}
(y-\langle y, x^*\rangle x) , x^*\right \rangle
= e^{it}\]
and
\[A(y)=i x +(i-\frac 1 2)(y-x).\]
So Theorem~\ref {A} does not hold if $X$ is not strictly convex.
\end{ex}

For any complex number 
$z \ne 0$, let $\mbox {sgn} \, (z)=\frac {\bar z}{|z|}$.

\begin{thm}  \label {C} Let $X$ be any strictly convex complex
Banach space and let $J:X \to X^*$ be any duality section of $X$.
If $x$ is a unit vector and if $T=\{e^{tA}:t \geq 0\}$ 
is a $(C_0)$ contraction
semigroup generated by $A$ such that
\[\lim_{t \to \infty} |\langle T(t) x, J(x)\rangle|=1,\]
then $x$ is an eigenvector of $A$ corresponding to a purely
imaginary eigenvalue.
\end{thm}

\begin{proof} We claim that for any $t$, $T(t) x= c(t) x$ for some
$|c(t)|=1$.

Suppose it is not true. Then there is $t \geq 0$ such that 
$f(d)=\|x+ d \, T(t) x\|$ is a continuous function from
$\{d \in \CC: |d|=1\}$  into the  interval $[0,2)$.
Hence, there is $\epsilon>0$ such that
\[ \|x + d \, T(t) x\| \leq 2 -\epsilon\]
for all $|d|=1$.  So if $s \geq 0$, then
\begin {eqnarray*}
2-\epsilon & \geq &  \left\|T(s)\left (x+
\mbox {sgn} \,\left (\frac {\langle
T(s+t)x, J(x) \rangle}{\langle T(s) x, J(x) \rangle}\right) T(t) x
\right)\right\| \\
& \geq &\left  |\left \langle T(s) x +\mbox {sgn} \, \left (\frac {\langle
T(s+t)x, J(x) \rangle}{\langle T(s) x, J(x) \rangle}\right) T(s+t) x
,J(x) \right \rangle \right |\\
&=&     \left |\langle T(s) x,J(x) \rangle\right|+\left|
\left\langle
\mbox {sgn} \,  \left(\frac {\langle
T(s+t)x, J(x) \rangle}{\langle T(s) x, J(x) \rangle}\right) T(s+t) x
,J(x) \right\rangle\right |.
\end{eqnarray*}
This contradicts that
\[\lim_{s \to \infty} |\langle T(s) x, J(x)\rangle|=1.\]
But   
\[A(x)=\lim_{t \to 0^+} \frac {e^{tA}(x)-x}
t=\lim_{t \to 0^+}\frac  {(c(t)-c(0))x} t  =c'(0) x.\]
This implies $x$ is also an eigenvector of $A$ corresponding to an
eigenvalue $c'(0)$.  
Hence,
\[1\geq\|e^A x\|= \|e^{c'(0)} x\|=|e^{c'(0)}|.\] 
On the other hand,
\[\|e^{A} x\| \geq \lim_{t \to \infty} \|T(t) x\| \geq \lim_{t
\to \infty} |\langle T(t) x,J(x)\rangle|=1.\]
This implies $|e^{c'(0)}|=1$, and so 
$c'(0)$ must be a purely imaginary number.  \end{proof}

Let $x$ be any unit vector in a complex Banach space such that
if $y \ne x$ and $\|y\|=1$, then $\|x+y\|<2$.  The proof above
shows that if $\{T=e^{tA}:t \geq 0\}$ is a $(C_0)$ contraction
semigroup generated by $A$ and if $\lim_{t \to \infty} |
\langle T(t)x ,
J(x)\rangle|=1$,
then $x$ is an eigenvector of $A$ corresponding to a purely
imaginary eigenvalue.

Let $x^*$ be any unit vector in $X^*$, and let $\alpha$ be any
real number between 0 and 1.  The {\em slice} $S(B,x^*, \alpha)$
is the set
\[ \{ x: \|x\|\leq 1, \mbox { and } \Re (\langle x, x^*\rangle)
\geq 1-\alpha\}.\]
Recall a unit vector $x$ is called an {\em exposed point} (of the
unit ball) if there is  $x^* \in X^*$ such that for any $y \ne
x$ with $\|y\| =1$,
\[1=\|x^*\|=\langle x,x^*\rangle>\Re \langle y, x^* \rangle.\]
The linear functional $x^*$ is said to {\em expose} $x$.  $x$ is called
a {\em strongly exposed  point} if there is  
$x^* \in X^*$ such that $\langle x,x^*\rangle=1$ and
\begin{eqnarray*}
0&=& \lim_{\alpha \to 0^+} \mbox {diam} \, S(B, x^*,\alpha)\\
  & =& \lim_{\alpha \to 0^+} \sup \{ \|y-z\|:y\in S(B,x^*,\alpha)
  \mbox { and } z \in S(B,x^*,\alpha)
   \}.\end{eqnarray*}
We say $x^*$  {\em strongly exposes} $x$.
One may ask the following question.

\begin{pro} Let $x$ be a unit vector.  Suppose that $x^*$ exposes
$x$.  
Is $x$  an eigenvector of $A$ corresponding to a purely
imaginary eigenvalue
if $T=\{e^{tA}:t \geq 0\}$ is a $(C_0)$ contraction
semigroup generated by $A$ such that
\[\lim_{t \to \infty} |\langle T(t) x, x^* \rangle|=1?\]
\end{pro}

We do not know the general answer to this question.  But the
following theorems show that the answer is affirmative 
if $x^*$ strongly exposes $x$
or if $X$ is a reflexive Banach space.

\begin{thm} \label {D} Let $x$ be a unit vector of a complex Banach space 
$X$.  Suppose that $x^*$
strongly exposes $x$.
If $T=\{e^{tA}:t \geq 0\}$ is a $(C_0)$ contraction
semigroup generated by $A$, and if
\[\lim_{t \to \infty} |\langle T(t) x, x^* \rangle|=1,\]
$x$ is an eigenvector of $A$ corresponding to a purely
imaginary eigenvalue.
\end{thm}

\begin{proof} By the proof of Theorem ~\ref {C}, it enough to
show that for each $s \geq 0$, $x$ is an eigenvector of $e^{sA}$.

Since $x^*$ strongly exposes $x$ and 
\[\lim_{t \to \infty} \mbox
{sgn}\,(\langle T(t) x, x^* \rangle) \langle T(t) x,
x^*\rangle=1,\]
$\mbox {sgn}\,(\langle T(t) x, x^* \rangle) T(t) x$ 
converges to $x$. Hence, for any $s \geq 0$,
\begin{eqnarray*}
e^{sA} x&=&e^{sA} (\lim_{t \to \infty} \mbox {sgn}\,
(\langle T(t) x, x^* \rangle) T(t) x)\\
&=&\lim_{t \to \infty}\mbox {sgn}\,
(\langle T(t) x, x^* \rangle)  T(s+t) x\\
&=&
\lim_{t \to \infty} \frac {\mbox {sgn}\,
(\langle T(t) x, x^* \rangle)}{\mbox {sgn}\,
(\langle T(s+t) x, x^* \rangle)}   x.\end{eqnarray*}
So $x$ is an eigenvector of $e^{sA}$.  We proved our theorem.
\end{proof}

\begin{thm}  Let $x$ be a unit vector of a complex reflexive Banach space 
$X$.  Suppose that $x^*$
exposes $x$.
If $T=\{e^{tA}:t \geq 0\}$ is a $(C_0)$ contraction
semigroup generated by $A$, and if
\[\lim_{t \to \infty} |\langle T(t) x, x^* \rangle|=1,\]
$x$ is an eigenvector of $A$ corresponding to a purely
imaginary eigenvalue.
\end{thm}

\begin{proof} 
As in the proof of Theorem ~\ref{D}, it is enough to show that $x$
is an eigenvector of $e^{sA}$ for all $s \geq 0$.

Fix $s \geq 0$ and a increasing positive sequence $\{t_n:n \in
\NN\}$ such that $\lim_{n \to \infty} t_n=\infty$.  Since $X$ is
reflexive, by passing to subsequences of $\{e^{t_n A} x: n \in
\NN\}$ and $\{e^{(s+t_n)A} x:n \in \NN\}$ we may assume that both
of them converge weakly, say
\begin{eqnarray*}
y& =&\mbox {w-}\lim_{n \to \infty} e^{t_n A} x  \\
z& = &\mbox {w-}\lim_{n \to \infty} e^{(s+t_n) A} x.
\end{eqnarray*}
But
\[\lim_{t \to \infty} |\langle e^{t A} x, x^*\rangle|=1,\]
and $x^*$ exposes $x$.  So if $c=\langle y, x^* \rangle
$,
and $d=\langle z,x^* \rangle $, then
$y=cx$, and
$z=dx$.  Hence
\begin{eqnarray*}e^{sA} (c x)&=& e^{sA} (\mbox
{w-} \lim_{n \to \infty} e^{t_n A} x)\\ &=& \mbox {w-}\lim_{n \to
\infty} e^{sA} ( e^{t_n A} x)\\ &=& dx \end{eqnarray*} 
This implies  $x$ is
an eigenvector of $A$ corresponding to a purely imaginary
eigenvalue. The proof is complete.\end{proof}

\begin{rem}  Let $x$ be a nonzero vector,
and let $T=\{T(t)=e^{tA}: t \geq 0\}$ be a $(C_0)$ contraction
semigroup generated by $A$ in a complex Banach space $X$.
Recently, Goldstein and Nagy \cite {GB} proved that if
$|\langle T(t) x, x^* \rangle| \to |\langle x, x^*\rangle| $ as
$t \to \infty$ for every $x^* \in X^*$, then $x$ is an
eigenvector of $A$ corresponding to a purely imaginary eigenvalue.
\end{rem}

\end{document}